\documentclass[12pt,letterpaper,titlepage]{amsart}
\usepackage{amsmath, amssymb, amsthm, amsfonts,amscd,xr}
\input epsf.tex

\newcommand{\R}{\mathbb{R}}
\newcommand{\C}{\mathbb{C}}

\def\beq{\begin{equation}}
\def\eeq{\end{equation}}

\def\arr{\hbox to 20pt{\rightarrowfill}}

\def\Gl{\mathrm {Gl}}

\def\L{\mathcal L}

\newenvironment{res} 
               {\begin{equation} 
\begin{minipage}{0.85\textwidth}} 
               { \end{minipage}\end{equation} } 
\def\ber{\begin{res} } 
\def\eer{\end{res}} 
 
\numberwithin{equation}{section} 
\newtheorem{thm}{Theorem}[section]

\newcommand{\oline}{\overline}

\newtheorem{lemma}[thm]{Lemma}

\newtheorem{prop}[thm]{Proposition}

\newtheorem{rem}[thm]{Remark}

\makeatletter 
\def\section{\@startsection {section}{1}{\z@}{3.5ex plus 1ex minus 
    .2ex}{2.3ex plus .2ex}{\large\bf}} 
    \def\subsection{\@startsection{subsection}{2}{\z@}{3.25ex plus 1ex minus 
 .2ex}{1.5ex plus .2ex}{\bf}} 
\makeatother 
\def\pf{{\em Proof}.\, } 
 
\def\bysame{\leavevmode\hbox to3em{\hrulefill}\,} 
 
\def\Ad{\operatorname{Ad}}

\def\af{\mathfrak{a}} 
\def\diag{\operatorname {diag}}
\def\gf{\mathfrak{g}} 
\def\h{\mathfrak{h}}

\def\im{\operatorname{im}}

\def\pf{\mathfrak{p}} 
\def\q{{\bf q}} 
\def\qf{\mathfrak{q}}

\def\L{\mathcal{L}}

\def\Sp{\mathrm{Sp}}
\def\PSp{\mathrm{PSp}}
\def\Sym{\mathrm{Sym}}
\def\PSl{\mathrm{PSl}}

\renewcommand{\Im}{\mbox{\rm Im}\,} 
\makeindex 
\begin{document} 
\title[Corner views] 
{Corner view on the crown domain}

\author{Bernhard Kr\"otz} 
\address{Max-Planck-Institut f\"ur Mathematik\\  
Vivatsgasse 7\\ D-53111 Bonn
\\email: kroetz@mpim-bonn.mpg.de}

\date{\today} 
\thanks{}
\maketitle 
 
\section {Introduction}

This paper is about the crown domain, henceforth denoted by $\Xi$, 
which is the
canonical complexification  
of a Riemmanian symmetric space $X$ of the non-compact type. 
Specifically we are interested in the the nature 
of the boundary  of $\Xi$ and eventually
in good compactifications of $\Xi$.

\par Let us begin with some possible definitions 
of $\Xi$. Let us denote by $G$ the connected component  
of the isometry group of $X$. Then  
$$X=G/K$$ 
for $K<G$ a maximal compact subgroup. Group-theoretically 
one can view $X$ as the moduli space of all maximal compact 
subgroups of the semisimple group  $G$. 

\par Before we advance let us recall some examples. 

\begin{itemize} 
\item For $G=\PSl(2,\R)$ it is custom to identify $X$ 
with the upper half plane
${\bf H}=\{z\in\C\mid \Im z>0\}$. Let us refer to 
\cite{K} for a comprehensive discussion of the 
corresponding crown domain and the analysis thereon.  
\item For $G=\PSl(n,\R)$  one can view $X$ as the space of  
positive definite unimodular  matrices, i.e.   
$X=\Sym(n,\R)_{\det =1}^+$.  
\item For $G=\PSp(n,\R)$ one often realizes $X$ as the Siegel 
domain $V+i\Omega$ with 
$V=\Sym(n,\R)$ and $\Omega\subset V$ the cone of positive definite 
matrices. 
\end{itemize}

\smallskip Next we discuss complexifications 
of $X$. By a complexification we simply mean a connected complex manifold 
$\Xi$ which contains $X$ as a totally real submanifold. 
We wish to request that complexification respects 
symmetry, i.e. the action of $G$ on $X$ extends 
to $\Xi$. 
\par The natural candidate for a complexification seems to 
be $X_\C=G_\C/K_\C$ with $G_\C$ and $K_\C$   
the universal complexifications of $G$ and $K$ respectively. We often call 
$X_\C$ the {\it affine complexification} of $X$ as it is an affine 
variety. For example for $X=\Sym(n,\R)_{\det =1}^+$ one has 
$X_\C=\Sym(n,\C)_{\det =1}$. 

\par Let us be more demanding on our complexification
and request that the Riemannian metric of $X$ extends 
to a $G$-invariant metric on $\Xi$ (see \cite{KSII}, 
Sect. 4). This is now a severe restriction on $\Xi$ 
as it forces the $G$-action on $\Xi$ to be proper. 
For instance $G$ does not act properly on $X_\C$
\footnote{The $G$-stabilizer of 
$\diag(i,-i)\in \Sym(2,\C)_{\det=1}$ is $\mathrm{PSO}(1,1)$ and not compact.}
and this guides us to  smaller $G$-domains 
in $X_\C$ where $G$ acts properly. 
The study of proper $G$-actions on $X_\C$ began in \cite{AG} and a complete 
classification of all maximal open $G$-neighborhoods of $X$  in $X_\C$ 
with proper action was obtained in \cite{K2}. 
In view of the results of \cite{K2} one can define the crown domain as 
the intersection of all maximal open 
$G$-neighborhoods of $X$  in $X_\C$ with proper action (\cite{K2}). 
\par Let us give another construction of $\Xi$. Denote 
by $TX$ the tangent bundle of $X$. This is a homogeneous
vector bundle over $X$ on which $G$ acts properly. There 
is a natural $G$-map from $TX$ to $X_\C$ (see \cite{KSII}, Sect. 4)
and the crown domain corresponds to  the maximal $G$-neighborhood of $X$ 
in $TX$ which embeds into $X_\C$. 

\par A third and perhaps prefered  way to define the crown domain is 
as the universal domain for holomorphically extended orbit maps 
of unitary  spherical representations of $G$ 
(see the introduction of \cite{KO}). 

\smallskip Let us turn now to the subject proper of this 
paper, the topological boundary $\partial\Xi$ of $\Xi$ in $X_\C$. 
The boundary is a very complicated object and there  is little hope 
to obtain an explicite description. However 
$\partial\Xi$ features some structure; for instance, in it one finds 
the distinguished boundary $\partial_d \Xi\subset\partial\Xi$, 
introduced in \cite{GK}. The distinguished boundary  is some sort of Shilov 
boundary of $\Xi$ in the sense that 
it is the smallest 
closed subset in $\partial\Xi$ on which bounded 
plurisubharmonic  functions on ${\rm cl}(\Xi)$ attend their maximum.

\par We know from \cite{GK} and \cite{KO} that 
$\partial_d\Xi$ is a finite (and explicite) union 
of $G$-orbits, say 
$$\partial_d\Xi= {\mathcal O}_1\amalg\ldots \amalg {\mathcal O}_s\, .$$
{}From  now on we shall identify each ${\mathcal O}_j$ with a homogeneous 
space: $G/H_j$. The main result of \cite{GK} was: 

\bigskip {\it If $G/H_j$ is a symmetric space, then it is 
a non-compactly causal symmetric space. Moreover, every non-compactly 
causal symmetric space $Y=G/H$ appears  in the distinguished boundary  
of the corresponding crown domain for $X=G/K$. }
\bigskip 

\par One  aim of this paper is to understand this 
result better. To be more concise: what is the reason that 
precisely non-compactly causal (NCC) symmetric spaces
appear in the boundary? As we will see, answering this question 
will eventually reveal the structure of $\partial\Xi$.

\par NCC-spaces are very special among all semisimple symmetric spaces. 
We recall their definition (see \cite{HO}). We assume 
the Lie algebra of $G$ to be simple and write $\qf$ for the tangent space 
of $Y$ at the standard base point $y_o=H\in Y$. 
We  note that 
$\qf$ is a linear $H$-module. Now, non-compactly causal means
that $\qf$ admits an non-empty open $H$-invariant convex cone, say $C$, 
which is hyperbolic and does not contain any affine lines.

\par The theme of this paper is to view 
$\Xi$ from the corner point $y_o\in Y$ and not as 
a thickening of $X$ as in the customary definitions from above.  
Now a slight precisioning  of terms is necessary.  
As we saw, $\partial_d\Xi$ might have several connected components. 
If this happens to be the case, then we shrink 
$\Xi$ to a $G$-domain $\Xi_H$ whose distinguished 
boundary is precisely $Y$, see \cite{GKO}. 

\par For $C\subset\qf$ the minimal cone (see \cite{HO}) 
we form in the tangent bundle 
$TY=G\times_H \qf$ the 
cone-subbundle 
$${\mathcal C}=G\times_H C $$ 
and with that its boundary cone-bundle 
$$\partial{\mathcal C} = G \times_H \partial C\, . $$  
In this context we  ask the following 

\bigskip{\it Question: Is there a $G$-equivariant, generically 
injective, proper continuous surjection $p: \partial {\mathcal C}
\to \partial \Xi_H$?}

\bigskip In other words,  we ask if there exists an equivariant 
"resolution" of the boundary in terms of the geometrically 
simple boundary cone bundle $\partial{\mathcal C}$. 

\par In this paper we give an affirmative answer 
to this question if $X$ is a Hermitian tube domain. 
In this simplified situation  
the crown domain is $\Xi=X\times \oline{X}$ with $\oline X$ 
denoting $X$ but endowed with the opposite complex structure
(i.e.,  if $X$ is already complex, then the crown 
is the complex double). On top of that $\partial_d\Xi=Y$ 
is connected, i.e. $\Xi=\Xi_H$. 

\par I wish to point out that the presented method of proof  
will not generalize. In order to advance one has to  
understand more about the structure of the minimal cone $C$; 
one might speculate that some sort of "$H\cap K$-invariant 
theory" for $C$ could be useful.   
\par Let me pose two open problems: 

\smallskip 

\noindent {\bf Problem 1:} For general $\Xi$, does 
$\partial\Xi_H$ admit a resolution
as a cone bundle in the sense described above. 
\smallskip
 
\noindent {\bf Problem 2:} Construct $G$-equivariant compactification
of $\Xi$, resp. $\Xi_H$.

\bigskip {\it Acknowledgement:} The origin of this paper 
traces back to my productive stay at the RIMS in 2005/2006. 
I am happy to express my gratitude to my former host 
Toshiyuki Kobayashi. Also I would like to thank 
Toshihiko Matsuki for some useful intuitive conversations 
arround this topic.

\section{Main part}

\par Let $X=G/K$ be a Hermitian 
symmetric space of tube type. This means that there is an
Euclidean (or formally real) Jordan algebra $V$ with positive cone 
$W\subset V $ such that 
$$X = V + iW \subset V_\C \, . $$
The action of $G$ is by fractional linear transformation 
and our choice of $K$ is such it fixes the base point 
$x_0= i e $ with $e\in V$ the identity element 
of the Jordan algebra. 

\par It is no loss of generality if we henceforth restrict ourselves 
to the basic case of $G=\Sp(n,\R)$ -- 
the more general case is obtained by using standard dictionary 
which can be found in text books, e.g. \cite{FK}.  

\par For our specific choice, the Jordan algebra 
is $V=\Sym(n,\R)$ and $W\subset V $ is the cone of positive 
definite symmetric matrices. The identity element $e$ 
is $I_n$, the  $n\times n$ identity matrix. 
The group $G$ acts on $X$ by standard fractional 
linear transformations:  $g=\begin{pmatrix} a & b \\ c & d\end{pmatrix}\in G $ 
with appropriate $a,\ldots, d\in M(n,\R)$ acts as 
$$ g\cdot z = (az +b)(cz+d)^{-1} \qquad (z\in X)\, .$$
The maximal compact subgroup $K$ identifies with $U(n)$ under 
the standard embedding 

$$U(n) \to G,\ \  u +i v \mapsto \begin{pmatrix} u & v \\ -v  & u \end{pmatrix}
\qquad (u,v  \in M(n,\R))\, .$$
It is  then clear that $K=U(n)$ is the stabilizer of $x_0= i I_n$. 
In the sequel we consider $V_\C$ as an  affine piece 
of the projective variety $\L$ of Lagrangians in $\C^{2n}$; 
the embedding is given by 
$$V_\C \mapsto \L, \ \ T\mapsto L_T:=\{ (T(v), v)\mid v\in \C^n\}\, .$$
It is then clear that $G_\C =\Sp(n,\C)$ acts on $\L$; in symbols: 
$g=\begin{pmatrix} a & b \\ c & d\end{pmatrix}\in G_\C $
 with appropriate $a,\ldots, d\in M(n,\C)$ acts as  

$$g\cdot L =\{ (av + bw, cv +dw)\mid (v,w)\in L\} \qquad(L\in \L).$$  
 
The space $\L$ is homogeneous under $G_\C$. If we 
choose the base point 
$$x_0 \leftrightarrow L_0=\{ (iv, v) \mid v \in \C^n\}, $$
then the stabilizer of $x_0$ in $G_\C$ is the  Siegel 
parabolic 
  
$$S^+= K_\C \ltimes P^+ \quad\hbox{and} \quad P^+=\left\{ {\bf 1}+ \begin{pmatrix} u & -iu \\ 
-i u & -u\end{pmatrix}\mid u \in V_\C\right\}\, .$$
Thus we have 
$$\L= G_\C \cdot L_0\simeq G_\C / S^+\, .$$ 
Sometimes it is useful to take the conjugate base point 
$\overline {x_0}= -i I_n$. Then the stabilizer of $\overline {L_0}$ in 
$\L$ is the oppposite Siegel parabolic

$$S^-= K_\C \ltimes P^- \quad\hbox{and} \quad P^-
=\left\{ {\bf 1}+ \begin{pmatrix} u & iu \\ 
i u & -u\end{pmatrix}\mid u \in V_\C\right\}$$
and     

$$\L= G_\C \cdot \overline{L_0}\simeq G_\C / S^-\, .$$ 

Next we come to the realization of the affine complexification 
$X_\C=G_\C/K_\C$. We consider the $G_\C$-equivariant 
embedding 

$$X_\C \to \L \times \L , \ \ gK_\C \mapsto (g\cdot L_0, 
g\cdot \overline{L_0})\, .$$ 
It is not hard to see that 

$$X_\C =\{ (L,L')\in \L\times \L\mid L+L'=\C^{2n}\}, $$
i.e., $X_\C$ is the affine variety of pairs of 
transversal Lagrangians. 

\par Set $\overline X = V-iW$ and note that the map 
$z\mapsto \overline z$ identifies 
$X$ with $\overline X$ in a $G$-equivariant, but antiholomorphic manner. 

\par Next we come to the subject matter, the {\it crown domain}
of $X$:
 
$$\Xi=X\times \oline X \subset X_\C\, .$$

Let us denote by $\partial\Xi$ the topological 
boundary of $\Xi$ in $X_\C$. The goal is to resolve $\partial \Xi$ by a cone bundle over 
the affine symmetric space $Y=G/H$ where $H=\Gl(n,\R)$ is the 
structure group of the Euclidean Jordan algebra $V$. 

\par We define an involution $\tau$ on $G$ by 

$$\tau(g)= I_{n,n} g I_{n,n}\quad \hbox{where}\quad 
I_{n,n}=\begin{pmatrix} I_n & \\ & -I_n\end{pmatrix}. $$
The fixed point set of $\tau$ is 
$$H=\left \{ \begin{pmatrix} a & \\ & a^{-t}\end{pmatrix}\mid 
a\in \Gl(n,\R)\right\} = \Gl(n,\R)\, .$$

We write $\h$ for the Lie algebra of $H$ and denote by $\tau$ as 
well the derived involution on $\gf$. 
The $\tau$-eigenspace decomposition on $\gf$ shall be denoted 
by 
$$\gf=\h+\qf\quad \hbox{where}\quad \qf=\begin{pmatrix} 0 & V \\ V & 0\end{pmatrix}\, .$$
Write $\qf^+=\begin{pmatrix} 0 & V \\ 0 & 0\end{pmatrix}$ and 
$\qf^-=\begin{pmatrix} 0 & 0 \\ V & 0\end{pmatrix}$ and note  that 
$$\qf=\qf^+\oplus \qf^-$$
is the splitting of $\qf$ into two inequivalent irreducible 
$H$-modules. 

\par The affine space $Y=G/H$ 
admits (up to sign) a unique $H$-invariant convex open cone $C\subset \qf$, containing 
no affine lines and consisting of hyperbolic elements. Explicitely: 

$$C=\begin{pmatrix} 0 & W \\ W & 0\end{pmatrix}= W \oplus W \subset
\q^+ \oplus \q^-\, .$$ 

We form the cone bundle 
$${\mathcal C}=G\times_H C $$
 and note that 
there is a natural $G$-equivariant map 

$$P: G\times_H C \to \Xi, \ \ [g, (y_1, y_2)]\mapsto 
g\cdot (iy_1, (iy_2)^{-1})\, .$$
Let us verify that this map is in fact defined. 
For that one needs to check that for $h\in H$ and $y_1, y_2\in W$, the elements 
$(h, y_1, y_2)$ and $({\bf 1}, hy_1h^t, h^{-t} y_2 h^{-1})$ have the 
same image. 
Indeed, 

$$h\cdot (iy_1, (iy_2)^{-1})=(i hy_1 h^t, h (iy_2)^{-1} h^t)
= (i hy_1 h^t,  (i h^{-t}y_2h^{-1})^{-1})$$
which was asserted. 

\begin{lemma}\label{l=1} The map $P: {\mathcal C} \to \Xi$ 
is onto. 
\end{lemma}

\begin{proof} Write $A$ for the group of diagonal matrices 
in $G$ with positive entries. Note that the Lie algebra
$\af$ of $A$ is a maximal flat in $\pf=\gf\cap \Sym(2n,\R)$. 
In general, we know that $\pf=\Ad(K) \af$. 
Furthermore, if $W_d$ denotes the diagonal part of $W$ , then 
$iW_d=A\cdot x_0$. From $G=KAK$ it now follows 
that for any two points $(z,w)\in X$ there exist a $g\in G$ 
such that $g\cdot (z,w)=(x_0, w')$ with $w'\in i W_d$. 
As a consequence we obtain that 
$$\Xi= G\cdot ( i W_d , -i I_n)\, .$$
Clearly the right hand side is contained in the image 
of $P$ and this finishes the proof. 
\end{proof}

\begin{rem} (a) The map $P$ is not injective. We shall give 
two different arguments for this assertion, beginning with 
an abstract one. 
If $P$ were injective, then $P$ establishes an homeoporphism 
between 
$\Xi$ and ${\mathcal C}= G\times _H C$. In particular $\Xi$ is homotopy equivalent 
to $Y=G/H$. But we know that $\Xi$ is contractible; a contradiction. 
\par More concretely for $k\in K, k\neq {\bf 1}$, the elements 
$[k, (iI_n, -iI_n)]\neq [{\bf 1},(iI_n, -i I_n)]$  have the same image 
in $\Xi$. 
It should be remarked however, that the map is generically injective. 
\par\noindent(b) As $H$ acts properly on $C$, it follows that 
$G$ acts properly on the cone-bundle $G\times_H C$. Further it is 
not hard to see that the map $P$ is proper. 
\end{rem}

We need a more invariant formulation of 
the map $P$. For that,  note that the rational map  
$$V_\C \to V_\C, \ \ z\mapsto - z^{-1}$$
belongs to $K$. Its extension to $\L$, shall be denoted 
by $s_0$ and  is given by 
$$ s_0(L)=\{ (-w, v)\in \C^{2n}\mid (v,w)\in L \}\,  .$$
Also, the  anti-symplectic map $V_\C\to V_\C,  z\mapsto - z$ has a
natural extension to $\L$ given by 
$$L\mapsto -L :=\{ (-v,w)\in \C^{2n}\mid (v,w)\in L\}\, .$$
In this way, we can rewrite $P$ as 
$$P: G\times_H C  \to \Xi, \ \ [g,(y_1,y_2)]\mapsto 
g\cdot(iy_1, -s_0(iy_2))$$
and we see that $P$ extends to a continuous map 
$$\tilde P: G\times_H \qf\to \L\times \L, \ \ [g,(y_1, y_2)]\mapsto 
 g\cdot (iy_1, - s_0(iy_2))\, .$$
We restrict $\tilde P$ to $G\times_H \partial C$ and 
call this restriction $p$. It is clear 
that $\im  p$ is contained in the 
boundary of $\Xi$ in $\L\times \L$. But even more is 
true: the following proposition 
constitutes a $G$-equivariant ``resolution'' of $\partial \Xi$.

\begin{prop} $\im p \subset \partial \Xi$ and the  $G$-equivariant map 
$$p: G\times_H \partial C \to \partial \Xi ,\ \  [g, (y_1,y_2)]\mapsto 
g\cdot (iy_1, -s_0(iy_2))$$
is onto and proper. 
\end{prop}

\begin{proof} We first show that $\im p \subset \partial \Xi$. 
This means that $\im p  \subset X_\C$. In fact, from Lemma \ref{l=1}
and the definition 
of $p$ it follows that $\im p $ is contained in the closure of $\Xi$ 
in $\L\times \L$ and does not intersect $\Xi$. 

\par Let us now show that $\im p  \subset X_\C$. First note 
that 
\begin{equation}\label{e1} \partial C = W \times \partial W \amalg 
\partial W \times \partial W \amalg \partial W \times W\, . 
\end{equation}
Thus the assertion will certainly follow if we verify 
the following slightly stronger statement: 
for $y_1,y_2\in cl(W) $ the Lagrangians
$$L_1=\{ (iy_1 v, v)\mid v\in \C^n\} \quad\hbox{and}\quad 
L_2=\{ (w, iy_2 w)\mid w\in \C^n\}$$
are transversal. We use the structure group $H$ to bring 
$y_1$ in normal form 
$$y_1=\mathrm{diag}(\underbrace{1,\ldots, 1}_{p-\hbox{times}}, 0,\ldots, 0)\,.$$
Thus $(iy_1 v, v) = (w, iy_2 w)$ for some $v,w\in\C^n$ means explicitely that  

$$(iv_1, iv_2, \ldots, iv_p, 0, \ldots, 0; v_1, \ldots, v_n)=
(w_1, \ldots, w_n; iy_2(w))\, .$$ 
We conclude that $w_{p+1}=\ldots= w_n=0$. If $p=0$, then we are finished. 
So let us assume that $p>0$. But then 
$$y_2=\begin{pmatrix} -I_p & * \\ * & *\end{pmatrix}$$
and this contradicts the fact that $y_2$ is positive semi-definite.

\par We turn our attention to the onto-ness of $p$. 
For that note that the closure $cl(X)$ in $\L$ equals the geodesic 
compactification. As a result $\partial X = K\cdot (i\partial W_d)
=K\cdot (i\partial W)$. Likewise 
$\partial \oline X= K\cdot (-i \partial W)$.  
Observe that 
\begin{equation} \label{e2}
\partial \Xi= \left[X \times \partial \oline X \amalg \partial X \times \partial 
\oline X \amalg \partial X \times \oline X\right]\cap X_\C\, .\end{equation}

We first show that $X\times \partial \oline X \subset \im p$, even more 
precisely $p(G\times_H (W \times \partial W ))=  X\times \partial \oline X$.
In fact, 
$$X\times \partial \oline X = G \cdot (iI_n, K \cdot i\partial W )
=G \cdot (iI_n,  i\partial W )$$
and the claim is implied by (\ref{e1}). In the manner one verifies 
that $\partial X \times \oline X \subset \im p$. 

\par In order to conclude the proof it is now enough to show 
that $p$ is proper. This is because proper maps 
are closed and we have already seen that $\im p$ contains 
the dense piece $X\times \partial \oline X \amalg \partial 
X\times \oline X \subset \partial \Xi$. 
Now to see that $p$ is proper, it is enough 
to show that inverse images of compact subsets 
in  $[\partial X\times \partial \oline X] \cap X_\C$ 
are compact. For the other pieces in $\partial \Xi$ this is more 
or less automatic: Use that $G$ acts properly on $X$, resp. $\oline X$
which implies that $G$ acts properly on  
$X\times \partial \oline X$ resp. $\partial X\times \oline X$; likewise $G$ acts properly 
on $G\times_H (W \times \partial W)$ and  
$G\times_H (\partial W \times W)$. 
Thus we are about to show that preimages of compacta in  
$[\partial X\times \partial \oline X] \cap X_\C$ are 
again compact. But this is more or less immediate 
from transversality; I allow myself to skip 
the details. 
\end{proof}

\begin{rem} For $n=1$ 
the map $p$ is in fact a homeomorphism which we showed in 
\cite{KO}.  If $n>1$, the map $p$ fails to be injective by the same 
computational reason shown in the preceeding remark. 
However, we emphasize that the map is generically injective and that 
$p|_{\partial C}$ is injective. 
\end{rem}

\end{document}